# The Hirsch function and its properties

Leo Egghe

Hasselt University, 3500 Hasselt, Belgium

leo.egghe@uhasselt.be

ORCID: 0000-0001-8419-2932

Abstract

The Hirsch function of a given continuous function is a new function depending on a parameter. It exists provided some assumptions are satisfied. If this parameter takes the value one, we obtain the well-known h-index. We prove some properties of the Hirsch function and characterize the shape of general functions that are Hirsch functions. We, moreover, present a formula that enables the calculation of f, given its Hirsch function $h_f$.

Keywords: h-index; h-function; Hirsch function

## 1. Introduction

Let $f: \mathbb{R}^+ \to \mathbb{R}^+$ be a function. Then we define for all $\theta \in \mathbb{R}_0^+ = \mathbb{R}^+ \setminus \{0\}$:

$$x = h_f(\theta) \Leftrightarrow f(x) = \theta x \qquad (1)$$

We only consider those cases for which (1) has a unique solution. If f(0) = 0 then we exclude a possible extra solution of $x = h_f(\theta) = 0$ unless this is a unique solution. Fig. 1 illustrates some special cases.

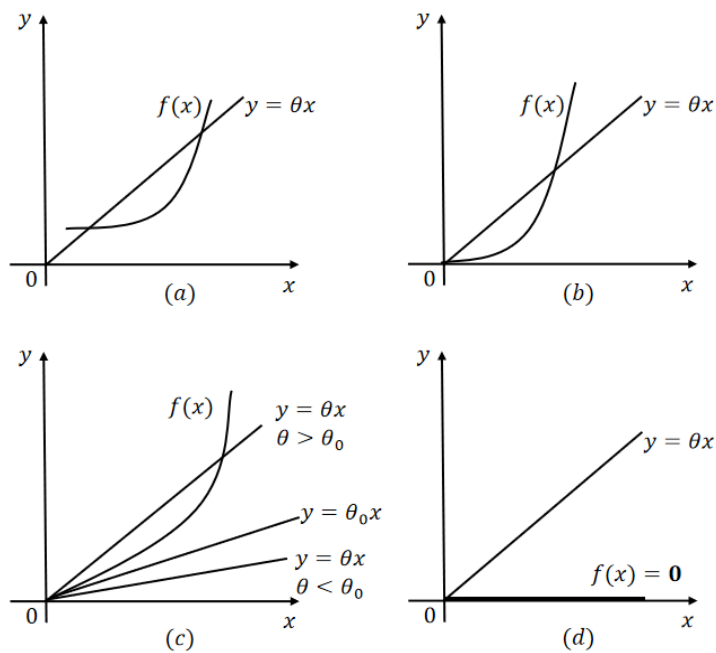

Fig.1. Some special cases

Case a) does not lead to a valid solution of (1) as y = θx and f(x) intersect in more than one point.

Case b). Here f(0) = f'(0) = 0. Here we do not consider x =0, so that (1) has a unique solution for all $\theta \in \mathbb{R}_0^+$.

Case c). Here f(0) = 0 and f'(0) = $\theta_0$ > 0. We do not consider x = 0 as a solution of (1) if θ > $\theta_0$ and do consider x = 0 as a solution if 0 < θ ≤ $\theta_0$.

Case d). Here we have x = $h_f(\theta)$ = 0, for all $\theta \in \mathbb{R}_0^+$.

Although it is possible to solve such special cases differently, the main point is that we know unambiguously what we mean by the notation $h_f(\theta)$. As $h_f(\theta)$ is now clearly defined we obtain a well-defined function $h_f$.

Definition: The Hirsch function

The function $h_f : \theta \in \mathbb{R}_0^+ \to h_f(\theta) \in \mathbb{R}^+$ is called the Hirsch function.

For θ =1 we obtain the well-known h-index (Hirsch, 2005) of the continuous function f, explaining the naming of this function. We further note that $h_f$ is not defined in point zero so we can say that for f = **0** (the null function) $h_f(\theta)$ = 0.

The Hirsch function has been used implicitly in (Egghe & Rousseau, 2019) (without naming it as such) and later in (Egghe, 2021, 2022), while the idea of considering h-indices with a variable parameter, originates from van Eck and Waltman (2008).

The Hirsch function is not defined as an explicit function but implicitly through equation (1). We first provide a characterization of such functions.

Theorem 1

Let φ be a function defined on $\mathbb{R}^+$, continuous in 0. Let further f be a function, continuous in the point φ(0) then the following two statements are equivalent:

(i) $h_f = \varphi$ on $\mathbb{R}_0^+$

(ii) $\forall \theta \in \mathbb{R}^+ : f(\varphi(\theta)) = \theta \cdot \varphi(\theta)$ (2)

Proof. (i) ⇒ (ii)

From (i) and (1) we obtain (2) $\forall \theta \in \mathbb{R}_0^+$. For θ = 0, we find, using the assumed continuity :

$$f(\varphi(0)) = f\left(\lim_{\theta \to 0} \varphi(\theta)\right) = \lim_{\theta \to 0}(f(\varphi(\theta))) = \lim_{\theta \to 0}(\theta \cdot \varphi(\theta)) = 0$$

Where we have used that we already know (2) for θ > 0. Hence, f(φ(0)) = 0 = 0.φ(0), which is (2) for θ = 0.

(ii) ⇒ (i)

From (2), (1), and the assumed uniqueness we have that $\forall \theta \in \mathbb{R}_0^+ : h_f(\theta) = \varphi(\theta)$, by the definition of $h_f$. □

Next, we will study the following problems

(a). Given f, determine $h_f$. This is the formalism shown in (1). One example: let f(x) = C > 0 (C fixed). Then (1) leads to the equation C = θx. Hence $\forall \theta \in \mathbb{R}_0^+$ $h_f(\theta) = x = C/\theta$. We come to the same result using (2). Indeed, $\forall \theta \in \mathbb{R}_0^+$ $h_f(\theta) = \varphi(\theta) = C/\theta$.

(b). Given φ, determine f such that φ = $h_f$. This problem already places some extra requirements on φ without which φ = $h_f$ is impossible. Consider e.g., the example above: with $\forall x \in \mathbb{R}_0^+$ φ(x) = C/x. Then (2) leads to $f\left(\frac{C}{x}\right) = x.\frac{C}{x} = C$. As the range of C/x is $\mathbb{R}_0^+$, f(x) = C on $\mathbb{R}_0^+$ and thus also f(x) = C on $\mathbb{R}^+$ by the continuity of f.

(c). Neither f nor φ is given, but a general relationship between f and φ. Here we consider two subcases.

1) φ is given via a relation with f

Example 1. φ = f (the simplest possible relation). Using (2) we have $\forall \theta \in \mathbb{R}^+$: f(f(θ)) = θ.f(θ). If f is continuous then this relation can only occur if f = **0** (the null function) or f(x) = $x^\alpha$, where α is the golden section, $\frac{\sqrt{5}+1}{2}$ (Egghe, Fibonacci article). Its proof uses the Fibonacci sequence.

Example 2. φ = f∘f. Using (2) this leads to f(f(f(θ))) = θ.f(f(θ)), $\forall \theta \in \mathbb{R}^+$, see (Egghe, Fibonacci). Again, for f continuous, this requirement leads to two possible solutions, namely f = **0** or f(x) = $x^\beta$, with β ≈ 1.4648493 (smaller than the golden section). Its proof uses a variant of the Fibonacci sequence.

In the same vein, one can consider the case $\varphi = \underbrace{f \circ f \circ \ldots \circ f}_{n \text{ times}}$.

2) f is given via a relationship with φ

Example 1. The function f = φ. Although this is essentially the same as the previous example 1, (2) leads to $\forall \theta \in \mathbb{R}^+$: φ(φ(θ)) = θ.φ(θ), leading to φ (=f) = **0** or φ(x) (=f(x)) = $x^\alpha$.

Example 2. f = φ∘φ

This example is different. Via (2) we find: φ(φ(φ(θ))) = θ.φ(φ(θ)), $\forall \theta \in \mathbb{R}^+$, see (Egghe, Fibonacci). For φ continuous this leads to φ = **0** or φ(x) = $x^\alpha$, with α ≈ 1.3247178. hence f(x) = $x^{(\alpha^2)}$.

This ends the introduction. Next, we will study the basic properties of the Hirsch function.

## 2. Properties of the Hirsch function

Theorem 2

The function $h_f$ is injective on the set $\{\theta \in \mathbb{R}_0^+ \,||\, h_f(\theta) \neq 0\}$.

Notation. We denote $\{\theta \in \mathbb{R}_0^+ \mid h_f(\theta) \neq 0\}$ as $\{h_f \neq 0\}$.

Proof. Let $x_1 = h_f(\theta_1) = h_f(\theta_2) = x_2$. Then (1) implies that $f(x_1) = \theta_1 x_1$ and $f(x_2) = \theta_2 x_2$. As $x_1 = x_2$ and f is a function this implies that $\theta_1 x_1 = \theta_2 x_2$, leading to $\theta_1 = \theta_2$ if $x_1 = x_2 \neq 0$. □

The next theorem provides a new characterization of $h_f$.

Theorem 3

Let m be a function of functions m: f → m(f), then the following statements are equivalent:

(i) $m(f) = h_f$

(ii) $\forall \theta \in \mathbb{R}_0^+ : m_\theta(f) = \psi_f^{-1}(\theta) = x$, where $\psi_f$ is injective, and defined as:

$$\psi_f(x) = \frac{f(x)}{x} = \theta \qquad (3)$$

Proof. (i) ⇒ (ii)

From (i) and (1) it follows that $\forall \theta \in \mathbb{R}_0^+ : x = m_\theta(f) \Leftrightarrow f(x) = \theta x \Leftrightarrow \theta = f(x)/x \Leftrightarrow x = \psi_f^{-1}(\theta)$. Moreover, from the fact that $h_f$ is a function, it follows that $\psi_f$ is injective.

(ii) ⇒ (i)

It follows similarly from (ii) and (1) that $m(f) = h_f$.

Remark

As $h_f^{-1} = \psi_f$, with $\psi_f$ defined in (3) it follows that $h_f^{-1}$ is a function on $\mathbb{R}_0^+$. This immediately leads to (see also Theorem 2):

$$\text{the function f is continuous} \Leftrightarrow h_f^{-1} \text{ is continuous} \qquad (4)$$

The two implications in (4) do not hold for $h_f$ (see further). To study this, we recall two results (stated as lemmas) from real analysis.

Lemma 1.

If f is continuous on an interval (possibly infinitely long) and injective then f is strictly monotonous.

Lemma 2

If f is injective, then the following two statements are equivalent:

(i) f is continuous on [a,b]

(ii) The function f⁻¹ is continuous on [f(a), f(b)] (or [f(b), f(a)])

A proof can be founded using Lemma 1 and (De Lillo, Theorem 2.27).

Notation. The domain of a function f is denoted as $\mathbf{D}(f)$.

Theorem 4

If $\mathbf{D}(f)$ is an interval, then f is continuous implies that $h_f$ is continuous.

Proof. $\mathbf{D}(\psi_f) = \mathbf{D}(f) \setminus \{0\}$, hence an interval. If f is continuous then also $\psi_f$ is continuous and $\psi_f$ is an injection (by Theorem 3). Applying now Lemma 2 on $\psi_f$ shows that $\psi_f^{-1}$ is a continuous function. It then follows from Theorem 3 that $h_f = \psi_f^{-1}$ is also continuous. □

Theorem 4 does not hold if one removes the requirement that $\mathbf{D}(f)$ is an interval. This is illustrated in Fig. 2.

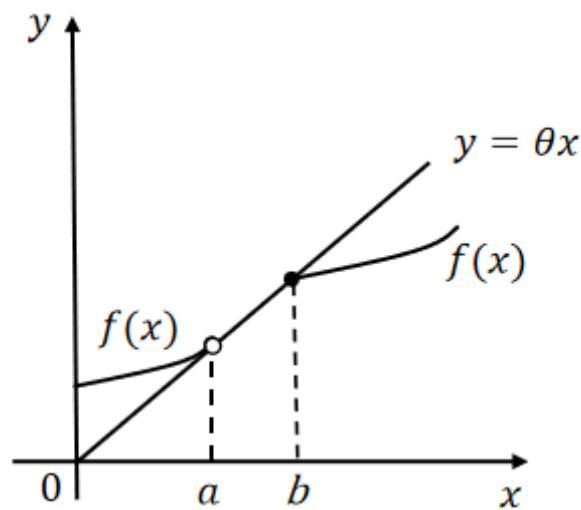

Fig. 2 A function f, continuous on its domain $\mathbf{D}(f)$ = [0,a[ ∪ [b, +∞ [ and a function $h_f(\theta)$ which is not continuous

We know that f is continuous if and only if $h_f^{-1}$ is continuous. Yet, we will show that the implication $h_f$ continuous ⇒ f continuous, does not always hold. For this, we need some preliminary results.

Lemma 3. If $f: \mathbb{R}^+ \to \mathbb{R}^+$ is continuous and injective on {f≠0}, then one of the following three statements hold:

(i) f is injective

(ii) ∃ $y_0 > 0$ such that $f|_{[0,y_0]} = \mathbf{0}$ and f (> 0) is injective, hence strictly increasing, on ]y₀, +∞ [,

(iii) ∃ $x_0 \geq 0$ such that $f|_{[x_0,+\infty[} = \mathbf{0}$ and f (> 0) is injective, hence strictly decreasing, on [0, x₀ [. Note that if x₀ = 0, this includes the case f = **0**.

Proof. Assume (i) is not the case, i.e., f is not injective. Yet, we know that f is injective on {f ≠ 0}. Hence, there exist x, y, 0 ≤ x < y such that f(x) = f(y) = 0.

We then show that f|[x,y] = **0**. (*)

Indeed, otherwise, there would exist z ∈ ]x,y[ such that f(z) ≠ 0. Because f is continuous it assumes all values between f(x) = 0 and f(z) > 0 on ]x,z[ and similarly on the interval ]z,y[. Consequently, there exist points x' and y', x'∈ ]x,z[ and y' ∈ ]z,y[, (hence x' ≠ y') such that f(x') = f(y') = f(z)/2 ≠ 0, which contradicts the fact that f is injective on {f ≠ 0}.

Next, we show that

$$\text{either } f|_{[0,x]} = \mathbf{0} \text{ or } f|_{[y,+\infty[} = \mathbf{0}. \quad (**)$$

Assume this is not the case. Then there exists u ∈ [0,x[ such that f(u) > 0 and v ∈ ]y, +∞[ such that f(v) > 0. As f is continuous it takes all values between f(x) = 0 en f(u) > 0 on ]u,x[ and between f(y) = 0 and f(v) > 0 on ]y,v[. Put a = min(f(u), f(v)) > 0. Then there exist x' in ]u,x[ and y' in ]y,v[ such that f(y') = f(x') = a ≠ 0 (and x' ≠ y'). This is in contradiction with the fact that f is injective on {f ≠ 0}.

From (*) and (**) it follows that f|[0,y] = **0** or f|[x, +∞[ = **0**, with 0 ≤ x < y. In the first case we set y₀ = sup{y > 0 such that f|[0,y] = **0**}. Then we know that $f|_{[0,y_0]} = \mathbf{0}$ and not equal to zero (hence strictly positive) on the compliment of [0,y₀]. From this, it follows that f is injective. Then f is strictly increasing on ]y₀, +∞[ and (ii) has been proved.

In the second case we set x₀ = inf{ x ≥ 0 such that f|[x,+∞[ = **0**}. Then $f|_{[x_0,+\infty[}$ = **0** and on the compliment of [x₀, +∞[ f ≠ 0 and hence injective. In this case, f decreases strictly on [0,x₀[ and (ii) is proved. □

Corollary

A continuous Hirsch function $h_f$ on $\mathbb{R}_0^+$ is of one of the following three types:

(i) $h_f$ is injective on $\mathbb{R}_0^+$;

(ii) $h_f = \mathbf{0}$ on an interval $[0, y_0]$, $y_0 > 0$ and strictly increasing on the compliment;

(iii) $h_f = \mathbf{0}$ on an interval $[x_0, +\infty[$, $x_0 > 0$ and strictly decreasing on the complement; including the case $h_f = \mathbf{0}$.

Proof. This follows from Lemma 3 and Theorem 2, with f (in Lemma 3) replaced by $h_f$, defined on $\mathbb{R}_0^+$. □

We show that these three types occur.

(i) This class is best known as it includes the functions $f(x) = x^c$, $c > 1$. Then, $\forall \theta \in \mathbb{R}_0^+ : h_f(\theta) = \theta^{\left(\frac{1}{(c-1)}\right)}$, which is a strictly increasing, injection. Note that for $1/(c-1) = c$ we find $c = \alpha$ and hence $f = h_f$ on $\mathbb{R}_0^+$.

(ii). See fig. 3a. The function f is strictly convex, $\theta_0 = f'(0) > 0$. By definition, $h_f$ is zero on $[0, \theta_0]$ and strictly increasing on $]\theta_0, +\infty[$.

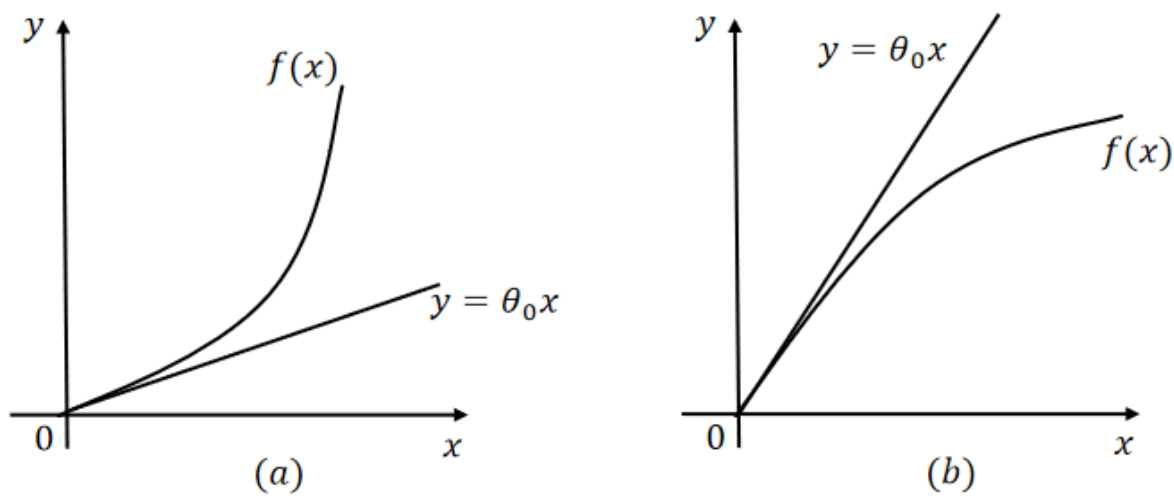

Fig. 3

(iii) See fig. 3b. The function f is strictly concave, $\theta_0 = f'(0) > 0$. By definition we have that $h_f$ is zero on $[\theta_0, +\infty[$ and strictly decreasing on $[0, \theta_0[$

The consequence of Lemma 3 also provides conditions for an equation such as (2), Theorem 1, to have or not to have a solution,

These cases are discussed in the next theorem.

Theorem 5.

If φ is a continuous function $\mathbb{R}^+ \to \mathbb{R}^+$ which is not of the form (i), (ii), or (iii) of the above corollary then a function f such that $h_f$ = φ does not exist. If φ is of the form (i), (ii), or (iii) then the solution of (2), namely $h_f$ = φ is given by

$$\forall x \in \mathbb{R}_0^+ : f(x) = x\, \varphi^{-1}(x) \qquad (5)$$

with $\varphi^{-1}$ the inverse function of the injective part of φ (abuse of notation). This function $\varphi^{-1}$ always exists, except when $x_0$ = 0 in (iii), in which case φ = **0** and φ = $h_f$ with f = **0**. We further note that in cases (ii) and (iii) f(0) = 0.

Proof.

Case (i). In this case, φ is injective and (2) gives:

$$\forall \theta \in \mathbb{R}_0^+ : f(\varphi(\theta)) = \theta\, \varphi(\theta)$$

Denoting φ(θ) by x we find that $\varphi^{-1}(x) = \theta$, which yields (5).

Case (ii). Now we know that there exists $y_0 > 0$ such that $\varphi|_{[0,y_0]}$ = **0** and φ is strictly increasing (hence injective) on $]y_0, +\infty[$. Next, we set $f(x) = x\, \varphi^{-1}(x)$ on $\varphi(]y_0, +\infty[)$. As φ is continuous and $\varphi|_{[0,y_0]}$ = **0,** $\mathbb{R}_0^+ \subset \varphi(]y_0, +\infty[)$. In this way, f is defined on $\mathbb{R}_0^+$ with (2) holding on $]y_0, +\infty[$. Now define f(0) = 0, then we have, $\forall \theta \in [0, y_0]$:

$$f(\varphi(\theta)) = f(0) = 0 = \theta\, \varphi(\theta)$$

showing that (2) holds on $\mathbb{R}^+$ and thus, by Theorem 1, $h_f$ = φ on $\mathbb{R}_0^+$.

Case (iii) Now we know that there exists $x_0 \geq 0$ such that $\varphi|_{[x_0,+\infty[}$ = **0** with φ strictly decreasing (and hence injective) on $[0, x_0[$. For $x_0$ = 0, φ = 0 on $\mathbb{R}^+$ and we take f = **0** on $\mathbb{R}^+$, leading to $h_f$ = φ on $\mathbb{R}_0^+$ (by (1)).

Assume now that $x_0 > 0$. Define $f(x) = x\, \varphi^{-1}(x)$ on $\varphi([0, x_0[) \neq \emptyset$. As φ is continuous and $\varphi|_{[x_0,+\infty[}$ = **0,** $\mathbb{R}_0^+ \subset \varphi([0, x_0[)$. So far, we defined f on $\mathbb{R}_0^+$, with (2) holding on $[0, x_0[$. Now, put f(0) = 0, then we have $\forall \theta \in [x_0, +\infty[$ :

$$f(\varphi(\theta)) = f(0) = 0 = \theta\, \varphi(\theta)$$

showing again that (2) holds on $\mathbb{R}^+$ and thus, by Theorem 1, $h_f$ = φ on $\mathbb{R}_0^+$. □

Practical conclusion

Leaving x = 0 aside we see that the solution of $h_f = \varphi$ is given by equation (5) with $\varphi^{-1}$ the inverse of $\varphi$ on the injective part of $\varphi$ (and f = **0** for $\varphi$ = **0**).

Examples

(i) For $\varphi(x) = C/x$, C>0 constant, we see that $\varphi$ is injective and $\varphi^{-1} = \varphi$. Then (5) yields: f(x) = x.C/x = C and $h_f = \varphi$.

For $\varphi(x) = x^c$, $\varphi^{-1}(x) = x^{1/c}$ and, by (5), $f(x) = x^{1+(\frac{1}{c})}$; $h_f = \varphi$.

(ii) and (iii). These cases are similar so we give just one example. For
$$\varphi(x) = \begin{cases} a^{x-b} - 1 \ (a > 1, b > 0), & x \geq b \\ 0 & 0 \leq x < b \end{cases}$$
we see that $\varphi$ is strictly increasing on [b, + ∞ [, and hence injective. On this set the function $\varphi^{-1}(x) = b + \log_a(x+1)$ and hence, using (5) we have:

$\forall x > 0$ ; f(x) = x (b + $\log_a$(x+1))  and f(0) = 0, showing that $h_f = \varphi$. Note that f'(0) = b and that $h_f$ is zero on [0,b].

Finally, we come to the case "$h_f$ continuous implies f continuous", the inverse statement of Theorem 4.

**Theorem 6**

If the range of f, denoted as **R**(f) is an interval, then $h_f$ continuous implies f continuous on **D**(f) ∩ $\mathbb{R}_0^+$ .

Proof

As $h_f$ is continuous everywhere, it is also continuous on $\{h_f \neq 0\}$, which is an interval  inside **D**($h_f$) = **D**($\psi_f^{-1}$), by the corollary to Lemma 3. By Theorem 3 **D**($\psi_f^{-1}$) = **R**($\psi_f$), which too is an interval because **R**(f) is an interval. By Theorem 2 we know that $h_f$ is injective on $\{h_f \neq 0\}$. Then it follows from Lemma 2 that $h_f^{-1} = \psi_f$ (by Theorem 3) is continuous on $\psi_f^{-1}(\{\psi_f^{-1} \neq 0\}) = $ **D**(f) ∩ $\mathbb{R}_0^+$. Finally, as $f(x) = x\, \psi_f(x)$ on $\mathbb{R}_0^+$ (by Theorem 3), this shows that f is continuous on **D**(f)  ∩ $\mathbb{R}_0^+$ .□

The next example shows that f is not necessarily continuous in zero. Take

$$f(x) = \begin{cases} x^2 & \text{for } x > 0 \\ 1 & \text{for } x = 0 \end{cases}$$

Then $\mathbf{R}(f) = \mathbb{R}_0^+$, which is an interval, $h_f$ is continuous on $\mathbb{R}_0^+$ but f is not continuous in 0, see Fig. 4.

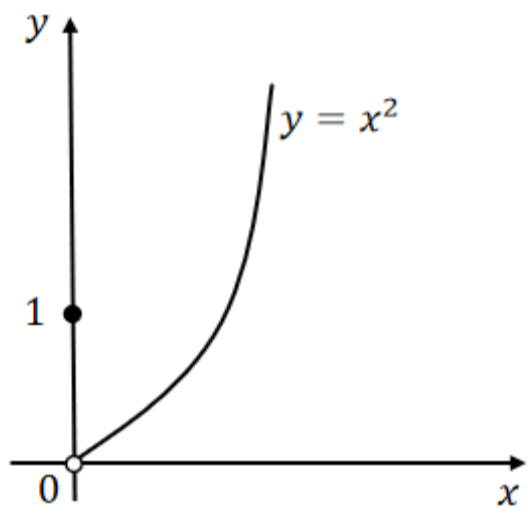

Fig. 4

We finish this article by remarking that the condition "$\mathbf{R}(f)$ is an interval" is necessary for Theorem 6. Fig. 5 provides an example of a function f which is not continuous on $\mathbb{R}_0^+$, but $h_f$ is continuous because $\mathbf{D}(h_f)$ is not an interval, (because $\mathbf{R}(f)$ is not an interval).

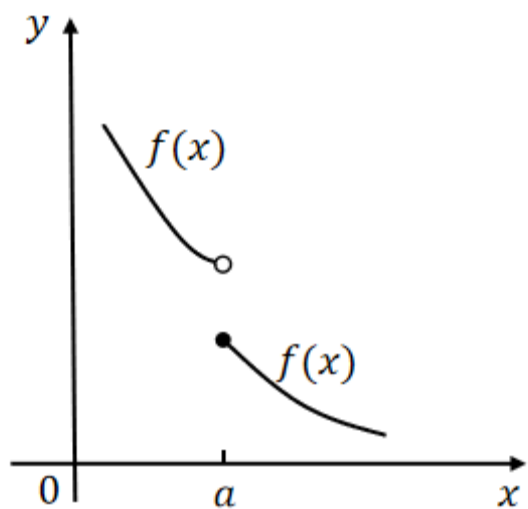

Fig. 5